\documentclass[10pt]{article}

\usepackage{graphics}
\usepackage{graphicx}

\usepackage[a4paper, left=35mm,right=35mm,top=34mm,bottom=34mm]{geometry}
\usepackage[utf8]{inputenc}
\usepackage[T1]{fontenc}
\usepackage[english]{babel}

\usepackage{enumerate}
\usepackage{graphicx}
\usepackage{hyperref}
\hypersetup{
    colorlinks=true,
    linkcolor=blue,
    filecolor=magenta,      
    urlcolor=cyan,
}

\usepackage{mathtools,amsthm,amssymb,amsfonts}
\usepackage{algorithm}
\usepackage{algorithmic}
\makeatother
\theoremstyle{plain}

\theoremstyle{definition}

\theoremstyle{remark}

\usepackage{caption} 
\usepackage[tight,footnotesize]{subfigure}

\usepackage{fancyhdr}
\usepackage{array}
\usepackage[usenames, dvipsnames]{xcolor}
\usepackage{tikz}
\usepackage{pgfplots}
\usepackage{ecp/ecpmashpc}
\usepackage{ecp/matrix_tikz/mrgcommon}
\usepackage{ecp/matrix_tikz/matrix_subdomain}

\author{
  {\normalsize Abal-Kassim Cheik Ahamed}\thanks{CUDA Research Center, Applied Mathematics and Systems Laboratory,  \'Ecole Centrale Paris, France.}
	\and
  {\normalsize Fr\'ed\'eric Magoul\`es}\thanks{CUDA Research Center, Applied Mathematics and Systems Laboratory, \'Ecole Centrale Paris, France
    (correspondence, frederic.magoules@hotmail.com).}
		}	
\title{Iterative Krylov Methods for Acoustic Problems on Graphics Processing Unit}
\date{}

\begin{document}
\maketitle
\thispagestyle{fancy}

\begin{abstract}
\noindent This paper deals with linear algebra operations on Graphics Processing Unit (GPU) with complex number arithmetic using double precision. An analysis of their uses within iterative Krylov methods is presented to solve acoustic problems.
Numerical experiments performed on a set of acoustic matrices arising from the modelisation of acoustic phenomena inside a car compartment are collected, and outline the performance, robustness and effectiveness of our algorithms, with a speed-up up to 28x for dot product, 9.8x for sparse matrix-vector product and solvers.
\end{abstract}


\section{Introduction}
\label{sec:introduction}
Linear algebra analysis has always been extremely useful when solving partial differential equations arising from many domains such as physics and biology models.
Even though Graphics Processing Units (GPUs) were first designed for graphic applications, they also represent a high potential for scientific computing and its applications to both physics and engineering. General-Purpose GPUs allow the developers to harness the high computational power of graphics cards to accelerate general-purpose scientific and engineering computing. The peak performance of CPUs and GPUs is significanlty different, due to the inherently different architectures between these processors.
In this work we focus on Compute Unified Device Architecture (CUDA)~\cite{GPU:CUDA4.0:2011}, proposed by NVIDIA in 2006, an appropriate and suitable language for NVIDIA graphics card. CUDA has offered a new vision in high performance computing.
In this paper, we analyse double precision complex number arithmetics algorithms of \emph{Alinea}~\cite{cheikahamed:2012:inproceedings-1,cheikahamed:2012:inproceedings-2}, our own research group library, which proposes linear algebra operations and iterative Krylov on both CPU and GPU clusters for real and complex number arithmetics in single and double precision.

The acoustic problem is steered in the frequency domain by the Helmholtz equation with suitable boundary conditions. The matrix of the linear system arising from the finite element discretization of the acoustic problem has a very huge size on high frequency regime.
Several discretization techniques like infinite element~\cite{magoules:journal-auth:26,magoules:journal-auth:19,magoules:journal-auth:15} or stabilized finite element~\cite{magoules:journal-auth:7} allows to reduce the size of the matrix.
The problem to solve comes from the discretization of the Helmholtz equation in a bounded domain $\Omega$, with outside boundary $\Gamma = \partial \Omega$. The Helmholtz equation is formulated as: $-\nabla^2 u - k^2 u= g$, where $k = \frac{2\pi F}{c}$ is the wavenumber of the frequency $F\in\mathbb{R}$ and $c\in\mathbb{R}$ is the velocity of the medium, which is different in space. In this work, we consider Dirichlet boundary conditions along a part of $\Gamma$.
Numerical experiments done on a set of acoustic finite element matrices are exhibited and show the performance, robustness and accuracy of linear algebra operations and their uses within iterative Krylov methods for solving acoustic problem modeled by Helmholtz equation.

The plan of this paper is the following. Section~\ref{sec:application_automotive_acoustic_flow} presents the industrial test cases involved to analyze our algorithms. Section~\ref{sec:linear_algebra_operations} presents numerical results of linear algebra operations required to carry out iterative Krylov methods such as addition of vectors, scale of vectors, sparse matrix-vector multiplication (SpMV), etc. Section~\ref{sec:iterative_krylov_methods} gives numerical tests on iterative Krylov methods, and Section~\ref{sec:concluding_remarks} gives conclusion.

\section{Application: automotive acoustic}
\label{sec:application_automotive_acoustic_flow}
\begin{figure}[!ht]
\centering
\includegraphics[scale=0.26]{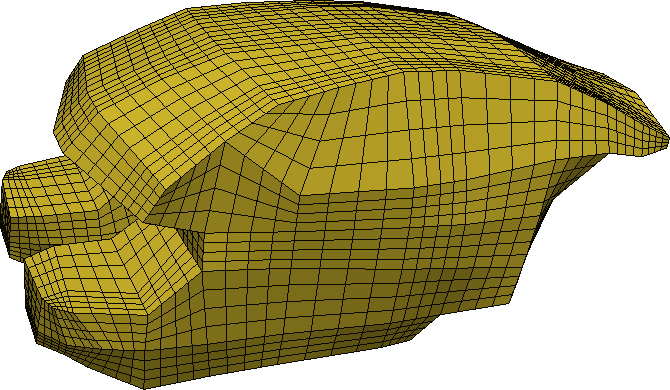}
\quad
\includegraphics[scale=0.26]{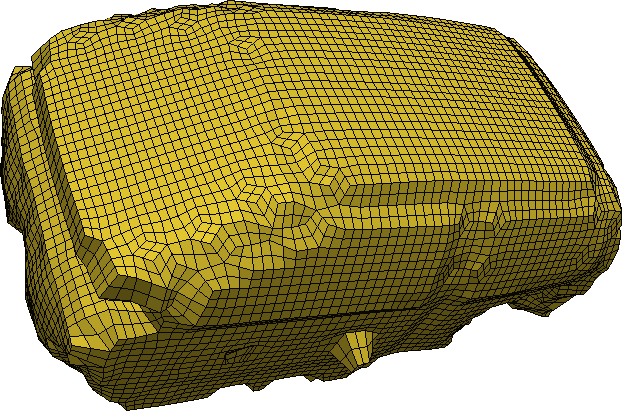}
\caption{Audi (Audi3D) and Twingo (Twingo3D)}
\label{fig:img:mesh}
\end{figure}

This part of the paper gives the main features of the finite element meshes used associated with the acoustic problems arising from the automotive industry~\cite{magoules:journal-auth:4}, namely car compartments: Audi (Audi3D) and Twingo (Twingo3D).
The car compartment problem is representative of acoustics cavity.
Fig.~\ref{fig:img:mesh} illustrates respectively the Audi3D and Twingo3D mesh for a given mesh size ($h$).

The matrices used to analyze and evaluate our algorithms are obtained from the finite element discretization of the acoustic problem, governed by the Helmholtz equation. The matrices are sparse large size, \emph{i.e.}, most values are zero. In this way, Compressed Sparse Row (CSR)~\cite{GPU:BG:2009}, is considered to store these matrices.
Table~\ref{tab:sketches_matrices_audi_Twingo_Cynlinder_3D} reports the matrices associated with the meshes of the car compartments. These features and characteristics are given in the third column. The sparse structure pattern and an histogramm of the distribution of nonzero values are respectively given in the first and second column.
\begin{table}
\caption{Sketches of Audi, Twingo FE matrices}
\label{tab:sketches_matrices_audi_Twingo_Cynlinder_3D}
\begin{tabular}{m{6.4cm}r}
\hlinewd{1.0pt}

\infosparse{./img/png/sparse/carHexOnlyFixed-GALERKIN-mat-1}{Audi3D-1}{1,727}{16,393}{0.550}{1,436}{27}{9.492}{10.205}\\
\multicolumn{2}{c}{\emph{3D acoustic FE matrix. Audi car (mesh size = 0.133425, length wave = 3.5).}}\\
\hlinewd{1.0pt}
\infosparse{./img/png/sparse/carHexOnlyFixed-GALERKIN-mat-2}{Audi3D-2}{11,637}{188,455}{0.139}{11,237}{27}{16.194}{11.223}\\
\multicolumn{2}{c}{\emph{3D acoustic FE matrix. Audi car (mesh size = 0.066604, length wave = 3.5).}}\\
\hlinewd{1.0pt}
\infosparse{./img/png/sparse/carHexOnlyFixed-GALERKIN-mat-3}{Audi3D-3}{85,001}{1,781,707}{0.025}{84,474}{27}{20.961}{9.832}\\
\multicolumn{2}{c}{\emph{3D acoustic FE matrix. Audi car (mesh size = 0.033289, length wave = 3.5).}}\\
\hlinewd{1.0pt}
\infosparse{./img/png/sparse/carHexOnlyFixed-GALERKIN-mat-4}{Audi3D-4}{648,849}{15,444,211}{0.004}{520,461}{27}{23.802}{7.720}\\
\multicolumn{2}{c}{\emph{3D acoustic FE matrix. Audi car (mesh size = 0.016643, length wave = 3.5).}}\\
\hlinewd{1.0pt}
\infosparse{./img/png/sparse/legacyCavity-GALERKIN-mat-0}{Twingo3D-0}{8,439}{143,889}{0.202}{6,268}{27}{17.050}{11.047}\\
\multicolumn{2}{c}{\emph{3D acoustic FE matrix. Twingo car (mesh size = 0.077866, length wave = 9.5).}}\\
\hlinewd{1.0pt}
\infosparse{./img/png/sparse/legacyCavity-GALERKIN-mat-1}{Twingo3D-1}{62,357}{1,351,521}{0.035}{53,935}{33}{21.674}{9.364}\\
\multicolumn{2}{c}{\emph{3D acoustic FE matrix. Twingo car (mesh size = 0.038791, length wave = 9.5).}}\\
\hlinewd{1.0pt}
\infosparse{./img/png/sparse/legacyCavity-GALERKIN-mat-2}{Twingo3D-2}{479,169}{11,616,477}{0.005}{470,625}{39}{24.243}{7.233}\\
\multicolumn{2}{c}{\emph{3D acoustic FE matrix. Twingo car (mesh size = 0.019379, length wave = 9.5).}}\\
\hlinewd{1.0pt}
\end{tabular}
\end{table}

The numerical experiments have been carried out on a workstation based on an Intel Core i7 920 2.67Ghz, which has 4 physical cores and 4 logical cores, 12GB RAM, and two NVIDIA graphics card: a Tesla K20c (device \#0) with 4799GB memory and GeForce GTX 570 with 1279MB memory (device \#1). The cards are double precision compatible. In the following Tesla K20c and GTX 570 will be denoted respectively \emph{gpu\#0} and \emph{gpu\#1}. For the sake of accuracy, we perform each operation 100 times, and the time indicated corresponds to the average time.

\section{Linear Algebra Operations}
\label{sec:linear_algebra_operations}

This section introduces linear algebra algorithms such as assign of a vector, scale of vectors, element wise product, addition of vectors, dot product and sparse matrix-vector products.
CUDA was originally dedicated for integer arithmetics and then for real numbers arithmetics, with a decreasing of performance of computations. Since, a complex number is a set of two real numbers composed of real and imaginary part, implementation is feasible by designing a structure of two real numbers. CUDA library includes a structure called \verb+cuComplex+, but for performance reasons, we specify our own complex class template structure \verb+complex<T>+ that offers all the operations given by the standard \verb+std::complex+.
As a result, in order to get the most benefits of GPU architecture, the elementary linear operation kernel requires to be reimplemented~\cite{journals:KnibbeOV11,GPU:BG:2009,GPU:BG:2008}.
In reference~\cite{cheikahamed:2012:inproceedings-2}, an analysis carried out on real number artihmetics with double precision with a suitable implementation of the CUDA kernel presents excellent speed-up for linear algebra operations and iterative Krylov methods~\cite{cheikahamed:2012:inproceedings-1}.
The finite element discretization of the Helmholtz equation for acoustic problems conducts to complex number arithmetics matrices.
In this paper, we give an extension of this analysis with acoustic problem. We develop efficient iterative Krylov methods for solving linear systems with complex number arithmetics.
As proved in~\cite{cheikahamed:2013:inproceedings-3} for real number arithmetics, our template code gives effective results compared to Cusp~\cite{GPU:CUSP:2010}, CUBLAS~\cite{GPU:CUBLAS}, CUSPARSE~\cite{GPU:CUSPARSE4.0:2011}.
But performance for complex number arithmetics with double precision remains a defiance, and dynamic auto-tuning of the GPU grid should be considered considered~\cite{cheikahamed:2012:inproceedings-2}.

The complex double precision running times in milliseconds (ms) of the \emph{assign operation} are collected in Table~\ref{tab:assign_tab}, with $h$ the size of the vector.
\begin{table}[!ht]  
\centering    
\renewcommand{\arraystretch}{1.3}
\renewcommand{\tabcolsep}{0.06cm}
\caption{Assign of vector (ZASSIGN)}
\label{tab:assign_tab}
\begin{tabular}{ccccccccc}
\hlinewd{1.0pt}
{\bf h} & {\bf cpu} & {\bf cpu} & {\bf gpu\#0} & {\bf gpu\#0} & {\bf gpu\#1} & {\bf gpu\#1} & {\bf ratio\#0} & {\bf ratio\#1} \\
{} & {\it time (ms)} & {\it Gflops} & {\it time (ms)} & {\it Gflops} & {\it time (ms)} & {\it Gflops} & {\it cpu/\#0} & {\it cpu/\#1} \\
\hlinewd{1.0pt}
648,849 & 1.10 & 0.59 & 0.16 & 4.06 & 0.20 & 3.27 & {\bf 6.88} & {\bf 5.54} \\
2,000,000 & 4.00 & 0.50 & 0.41 & 4.92 & 0.46 & 4.36 & {\bf 9.84} & {\bf 8.72} \\
9,000,000 & 18.33 & 0.49 & 1.79 & 5.04 & 1.69 & 5.31 & {\bf 10.27} & {\bf 10.82} \\
14,000,000 & 27.50 & 0.51 & 2.63 & 5.32 & 2.86 & 4.90 & {\bf 10.45} & {\bf 9.63} \\
\hlinewd{1.0pt}
\end{tabular}
\end{table}

In the following, all kernels compute the global index of each thread as follows:
\begin{lstlisting}[language=cuda,label=lst:global_idx_2d]
  unsigned int x = blockIdx.x * blockDim.x + threadIdx.x;
  unsigned int y = threadIdx.y + blockIdx.y * blockDim.y;
  int pitch = blockDim.x * gridDim.x;
  int idx = x + y * pitch;
\end{lstlisting}
The scale {\em scale operation} kernel is described as follows
\begin{lstlisting}[language=cuda,label=lst:xmy_2d]
__global__ void Scal( stdmrg::complex<double> alpha,
                      const stdmrg::complex<double>* d_x, int size) {
  if ( idx < size ) d_x[idx] = alpha * d_x[idx];
}
\end{lstlisting}
In Table~\ref{tab:scal_vectors_tab}, we collect the execution times of the scale {\em scale operation}.
\begin{table}[!ht]  
\centering    
\renewcommand{\arraystretch}{1.3}
\renewcommand{\tabcolsep}{0.06cm}
\caption{Scale of vectors (ZSCAL)}
\label{tab:scal_vectors_tab}
\begin{tabular}{ccccccccc}
\hlinewd{1.0pt}
{\bf h} & {\bf cpu} & {\bf cpu} & {\bf gpu\#0} & {\bf gpu\#0} & {\bf gpu\#1} & {\bf gpu\#1} & {\bf ratio\#0} & {\bf ratio\#1} \\
{} & {\it time (ms)} & {\it Gflops} & {\it time (ms)} & {\it Gflops} & {\it time (ms)} & {\it Gflops} & {\it cpu/\#0} & {\it cpu/\#1} \\
\hlinewd{1.0pt}
648,849 & 5.56 & 0.70 & 0.20 & 19.47 & 0.21 & 18.34 & {\bf 27.78} & {\bf 26.17} \\
2,000,000 & 15.71 & 0.76 & 0.46 & 26.04 & 0.53 & 22.56 & {\bf 34.10} & {\bf 29.54} \\
9,000,000 & 80.00 & 0.68 & 1.92 & 28.08 & 2.33 & 23.22 & {\bf 41.60} & {\bf 34.40} \\
14,000,000 & 120.00 & 0.70 & 2.94 & 28.56 & 3.57 & 23.52 & {\bf 40.80} & {\bf 33.60} \\
\hlinewd{1.0pt}
\end{tabular}
\end{table}

\noindent {\em Double-precision complex Alpha X Plus Y (Zaxpy)}, i.e., $y[i] = \alpha \times x[i] + y[i]$, is a level one (vector) operation between two complex number arithmetics vectors in the Basic Linear Algebra Subprograms (BLAS) package. The simple CUDA kernel of Zaxpy is implemented as follows:
\begin{lstlisting}[language=cuda,label=lst:zaxpy_2d]
__global__ void Daxpy(stdmrg::complex<double> alpha,
                      const stdmrg::complex<double>* d_x,
                      stdmrg::complex<double>* d_y, int size) {
  if ( idx < size ) d_y[idx] = alpha * d_x[idx] + d_y[idx];
}
\end{lstlisting}
In Table~\ref{tab:addition_vectors_tab}, we present the complex number arithmetics with double precision execution times in milliseconds (ms) of \emph{Zaxpy} operation.
\begin{table}[!ht]  
\centering    
\renewcommand{\arraystretch}{1.3}
\renewcommand{\tabcolsep}{0.06cm}
\caption{Addition of vectors (ZAXPY)}
\label{tab:addition_vectors_tab}
\begin{tabular}{ccccccccc}
\hlinewd{1.0pt}
{\bf h} & {\bf cpu} & {\bf cpu} & {\bf gpu\#0} & {\bf gpu\#0} & {\bf gpu\#1} & {\bf gpu\#1} & {\bf ratio\#0} & {\bf ratio\#1} \\
{} & {\it time (ms)} & {\it Gflops} & {\it time (ms)} & {\it Gflops} & {\it time (ms)} & {\it Gflops} & {\it cpu/\#0} & {\it cpu/\#1} \\
\hlinewd{1.0pt}
648,849 & 5.56 & 0.93 & 0.26 & 20.04 & 0.27 & 19.52 & {\bf 21.44} & {\bf 20.89} \\
2,000,000 & 16.67 & 0.96 & 0.69 & 23.20 & 0.81 & 19.68 & {\bf 24.17} & {\bf 20.50} \\
9,000,000 & 75.00 & 0.96 & 3.03 & 23.76 & 3.33 & 21.60 & {\bf 24.75} & {\bf 22.50} \\
14,000,000 & 120.00 & 0.93 & 4.76 & 23.52 & 5.26 & 21.28 & {\bf 25.20} & {\bf 22.80} \\
\hlinewd{1.0pt}
\end{tabular}
\end{table}

\noindent The {\em element wise product} or {\em element by element product}, i.e., $y[i] = x[i] \times y[i]$. The CUDA kernel, is described simply as: 
\begin{lstlisting}[language=cuda,label=lst:xmy_2d]
__global__ void EWProduct( stdmrg::complex<double> alpha,
  const stdmrg::complex<double>* d_x,stdmrg::complex<double>* d_y, int size) {
  int idx = x + y * pitch;
  if ( idx < size ) d_y[idx] = d_x[idx] * d_y[idx];
}
\end{lstlisting}
Table~\ref{tab:ewp_vectors_tab} exhibits the double precision execution times of the \emph{element by element product} operation.
\begin{table}[!ht]  
\centering    
\renewcommand{\arraystretch}{1.3}
\renewcommand{\tabcolsep}{0.10cm}
\caption{Element wise product (ZAXMY)}
\label{tab:ewp_vectors_tab}
\begin{tabular}{ccccccccc}
\hlinewd{1.0pt}
{\bf h} & {\bf cpu} & {\bf cpu} & {\bf gpu\#0} & {\bf gpu\#0} & {\bf gpu\#1} & {\bf gpu\#1} & {\bf ratio\#0} & {\bf ratio\#1} \\
{} & {\it time (ms)} & {\it Gflops} & {\it time (ms)} & {\it Gflops} & {\it time (ms)} & {\it Gflops} & {\it cpu/\#0} & {\it cpu/\#1} \\
\hlinewd{1.0pt}
648,849 & 8.33 & 0.47 & 0.28 & 13.66 & 0.29 & 13.55 & {\bf 29.25} & {\bf 29.00} \\
2,000,000 & 25.00 & 0.48 & 0.72 & 16.56 & 0.85 & 14.16 & {\bf 34.50} & {\bf 29.50} \\
9,000,000 & 120.00 & 0.45 & 3.03 & 17.82 & 3.33 & 16.20 & {\bf 39.60} & {\bf 36.00} \\
14,000,000 & 180.00 & 0.47 & 4.76 & 17.64 & 5.00 & 16.80 & {\bf 37.80} & {\bf 36.00} \\
\hlinewd{1.0pt}
\end{tabular}
\end{table}

\noindent Dot product operation can be very costly for large size vectors. 
Instead of performing a simple loop with simultaneous sums to compute the dot product, which is not very effective on GPUs, we perform it into two distinct tasks. The first is the element wise product of vectors and the second consists in summing all the results obtained at the first step. The reduction done at the second step associates each element of the input data with a thread, and at the end the partial sum of the $n^{th}$ first elements is stored in the first thread of the current block. The final dot product result is then computed as the sum of all the partial sums of the different blocks.
The double precision execution times of the \emph{dot product} on both CPU and GPU are exposed in Table~\ref{tab:dot_product_tab} and Fig.~\ref{fig:tikz:z_dot}.
\begin{table}[!ht]  
\centering    
\renewcommand{\arraystretch}{1.3}
\renewcommand{\tabcolsep}{0.06cm}
\caption{Dot product (ZDOT)}
\label{tab:dot_product_tab}
\begin{tabular}{ccccccccc}
\hlinewd{1.0pt}
{\bf h} & {\bf cpu} & {\bf cpu} & {\bf gpu\#0} & {\bf gpu\#0} & {\bf gpu\#1} & {\bf gpu\#1} & {\bf ratio\#0} & {\bf ratio\#1} \\
{} & {\it time (ms)} & {\it Gflops} & {\it time (ms)} & {\it Gflops} & {\it time (ms)} & {\it Gflops} & {\it cpu/\#0} & {\it cpu/\#1} \\
\hlinewd{1.0pt}
648,849 & 5.56 & 0.93 & 0.33 & 15.83 & 0.33 & 15.94 & {\bf 16.94} & {\bf 17.06} \\
2,000,000 & 16.67 & 0.96 & 0.83 & 19.20 & 0.76 & 20.96 & {\bf 20.00} & {\bf 21.83} \\
9,000,000 & 80.00 & 0.90 & 3.23 & 22.32 & 3.23 & 22.32 & {\bf 24.80} & {\bf 24.80} \\
14,000,000 & 130.00 & 0.86 & 4.76 & 23.52 & 4.55 & 24.64 & {\bf 27.30} & {\bf 28.60} \\
\hlinewd{1.0pt}
\end{tabular}
\end{table}
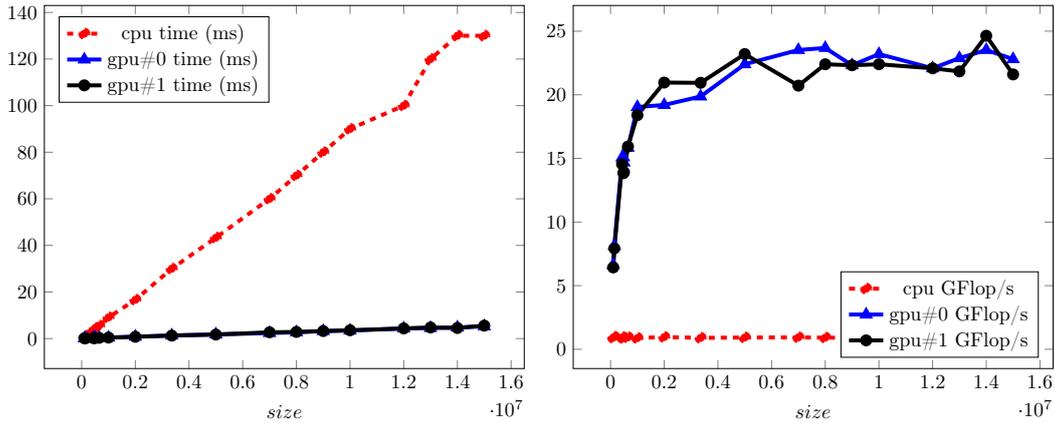
\begin{figure}[!ht]
\centering
  \begin{tikzpicture}[scale=0.75]
    \begin{axis}[
      height=8cm,
      width=10cm,
      xlabel=$size$,
      legend pos=north west,
      enlargelimits,
    ]
    \addplot[line width=2pt,dashed,red,mark=*] table[x index=0,y index=1,col sep=comma] {data/tikz/z_dot.txt};
    \addlegendentry{cpu time (ms)}
    \addplot[line width=2pt,blue,mark=triangle] table[x index=0,y index=3,col sep=comma] {data/tikz/z_dot.txt};
    \addlegendentry{gpu\#0 time (ms)}
    \addplot[line width=2pt,black,mark=*] table[x index=0,y index=5,col sep=comma] {data/tikz/z_dot.txt};
    \addlegendentry{gpu\#1 time (ms)}
    \end{axis}
  \end{tikzpicture}
  \begin{tikzpicture}[scale=0.75]
    \begin{axis}[
      height=8cm,
      width=10cm,
      xlabel=$size$,
      legend pos=south east,
      enlargelimits,
    ]
    \addplot[line width=2pt,dashed,red,mark=*] table[x index=0,y index=2,col sep=comma] {data/tikz/z_dot.txt};
    \addlegendentry{cpu GFlop/s}
    \addplot[line width=2pt,blue,mark=triangle] table[x index=0,y index=4,col sep=comma] {data/tikz/z_dot.txt};
    \addlegendentry{gpu\#0 GFlop/s}
    \addplot[line width=2pt,black,mark=*] table[x index=0,y index=6,col sep=comma] {data/tikz/z_dot.txt};
    \addlegendentry{gpu\#1 GFlop/s}
    \end{axis}
  \end{tikzpicture}
\caption{ZDOT [left: time in ms, right: GFlops]}
\label{fig:tikz:z_dot}
\end{figure}
Table~\ref{tab:normL2_tab} gives the numerical results of the norm operation.
\begin{table}[!ht]  
\centering    
\renewcommand{\arraystretch}{1.3}
\renewcommand{\tabcolsep}{0.06cm}
\caption{NormL2 (ZNORM)}
\label{tab:normL2_tab}
\begin{tabular}{ccccccccc}
\hlinewd{1.0pt}
{\bf h} & {\bf cpu} & {\bf cpu} & {\bf gpu\#0} & {\bf gpu\#0} & {\bf gpu\#1} & {\bf gpu\#1} & {\bf ratio\#0} & {\bf ratio\#1} \\
{} & {\it time (ms)} & {\it Gflops} & {\it time (ms)} & {\it Gflops} & {\it time (ms)} & {\it Gflops} & {\it cpu/\#0} & {\it cpu/\#1} \\
\hlinewd{1.0pt}
648,849 & 11.11 & 0.29 & 0.31 & 10.54 & 0.26 & 12.65 & {\bf 36.11} & {\bf 43.33} \\
2,000,000 & 33.33 & 0.30 & 0.73 & 13.70 & 0.57 & 17.60 & {\bf 45.67} & {\bf 58.67} \\
9,000,000 & 150.00 & 0.30 & 3.13 & 14.40 & 2.27 & 19.80 & {\bf 48.00} & {\bf 66.00} \\
14,000,000 & 230.00 & 0.30 & 5.00 & 14.00 & 3.70 & 18.90 & {\bf 46.00} & {\bf 62.10} \\
\hlinewd{1.0pt}
\end{tabular}
\end{table}

As shown in Table~\ref{tab:dot_product_tab} and Table~\ref{tab:normL2_tab}, GPUs clearly show better results than CPU with complex number arithmetics in double precision.
Much more than the dot product, the SpMV is probably the most time consuming operation in sparse matrix computation. This is required on all Krylov iterative methods. As proved in~\cite{cheikahamed:2012:inproceedings-2}, the performance of SpMV strongly depends on the properties of the matrix, particularly on the distribution of nonzero values. The following results are obtained with advanced auto-tuned techniques to organize threads on the CUDA grid.
References~\cite{GPU:KK:2011,GPU:ZR:2012,GPU:ZR:2010,GPU:KGWHBA:2011,GPU:DFG:2010} clearly showed the effectiveness of SpMV on GPU compared to CPU for real number arithmetics.
The running time and the number of floating operations per second for SpMV with complex number arithmetics with double precision are reported in Table~\ref{tab:csr_spmv_tab}.
\begin{table}[!ht]  
\centering    
\renewcommand{\arraystretch}{1.3}
\renewcommand{\tabcolsep}{0.07cm}
\caption{SpMV CSR}
\label{tab:csr_spmv_tab}
\begin{tabular}{lcccccccc}
\hlinewd{1.0pt}
{\bf problem} & {\bf cpu} & {\bf cpu} & {\bf gpu\#0} & {\bf gpu\#0} & {\bf gpu\#1} & {\bf gpu\#1} & {\bf ratio\#0} & {\bf ratio\#1} \\
{} & {\it time (ms)} & {\it Gflops} & {\it time (ms)} & {\it Gflops} & {\it time (ms)} & {\it Gflops} & {\it cpu/\#0} & {\it cpu/\#1} \\
\hlinewd{1.0pt}
Audi3D-0 & 0.01 & 0.61 & 0.07 & 0.12 & 0.06 & 0.13 & {\bf 0.19} & {\bf 0.21} \\
Audi3D-1 & 0.20 & 0.67 & 0.11 & 1.23 & 0.12 & 1.07 & {\bf 1.84} & {\bf 1.60} \\
Audi3D-2 & 2.22 & 0.68 & 0.37 & 4.03 & 0.42 & 3.56 & {\bf 5.93} & {\bf 5.24} \\
Audi3D-3 & 20.00 & 0.71 & 2.22 & 6.41 & 3.03 & 4.70 & {\bf 9.00} & {\bf 6.60} \\
Audi3D-4 & 180.00 & 0.69 & 18.33 & 6.74 & 24.00 & 5.15 & {\bf 9.82} & {\bf 7.50} \\
Twingo3D-0 & 1.67 & 0.69 & 0.28 & 4.06 & 0.33 & 3.45 & {\bf 5.88} & {\bf 5.00} \\
Twingo3D-1 & 15.71 & 0.69 & 1.79 & 6.05 & 2.44 & 4.43 & {\bf 8.80} & {\bf 6.44} \\
Twingo3D-2 & 140.00 & 0.66 & 14.29 & 6.51 & 16.67 & 5.58 & {\bf 9.80} & {\bf 8.40} \\
\hlinewd{1.0pt}
\end{tabular}
\end{table}

\section{Iterative Krylov methods}
\label{sec:iterative_krylov_methods}
After the analysis of linear algebra operations for complex number arithmetics with double precision, we now evaluate and analyze their uses within iterative Krylov methods~\cite{GPU:LS:2010,GPU:KK:2011,GPU:BCK:2011,GPU:ZR:2012}.
We have thus implemented a preconditionned bi-conjugate gradient stabilized method (P-Bi-CGSTAB), a preconditionned P-BiCGSTAB parametered (l) and a preconditionned transpose-free quasi-minimal residual method (P-tfQMR)~\cite{saad_iterative_2003}, with optimized CUDA and dynamic auto-tuning on GPU.
The data transfer between CPU and GPU consists of an important part of optimization~\cite{GPU:CSVGM:2014} for optimal performance on GPGPU. In our Krylov methods codes, we take care to send once all required input data from CPU to GPU before beginning the iterations. Even so, at each computed dot product or norm, there is one copy back from GPU to CPU. Both CPU and GPU codes are strictly the same, but all linear algebra operations such as Zdot, Znorm, Zaxpy, or SpMV are performed on device (GPU) for the GPU version.
The presented iterative Krylov methods are performed with a residual tolerance threshold of $1\times10^{-9}$, an initial guess of zero and $1000$ maximum number of iterations.
The numerical experiments presented in the following give an analysis of Krylov methods on CPU and GPU, with the same code, for complex number arithmetics with double precision.
The CPU and GPU execution times and corresponding speed-up of Audi3D and Twingo3D are collected in Table~\ref{speed-up_Audi3D} and Table~\ref{speed-up_Twingo3D}.
\begin{table}[htbp]
\setlength\columnsep{0.1pt}
\centering
\caption{Speed-up of Audi3D}
\label{speed-up_Audi3D}
\begin{tabular}{lcccc}
\hlinewd{1.0pt}
problem & \#iter & CPU time (s) & GPU time (s) & speed-up \\
\hline
\multicolumn{3}{l}{\emph{P-BiCGSTAB}} & & \\
Audi3D-1 & 21 & 0.01 & 0.030 & \textbf{0.33} \\
Audi3D-2 & 53 &  0.24 & 0.106 & \textbf{2.26} \\
Audi3D-3 & 94 &  4.01 & 0.703 & \textbf{5.71} \\
Audi3D-4 & 183 & 85.70 & 9.209 & \textbf{9.31} \\
\hline
\multicolumn{3}{l}{\emph{P-BiCGSTAB(8)}} & & \\
Audi3D-1 & 6 & 0.03 & 0.110 & \textbf{0.27} \\
Audi3D-2 & 12 & 0.52 & 0.286 & \textbf{1.82} \\
Audi3D-3 & 31 & 12.47 & 2.162 & \textbf{5.77} \\
Audi3D-4 & 70 & 266.26 & 30.100 & \textbf{8.85} \\
\hline
\multicolumn{3}{l}{\emph{P-TFQMR}} & & \\
Audi3D-1 & 24 & 0.02 & 0.040 & \textbf{0.50} \\
Audi3D-2 & 52 & 0.27 & 0.113 & \textbf{2.40} \\
Audi3D-3 & 99 &  4.71 & 0.755 & \textbf{6.24} \\
Audi3D-4 & 214 &  102.17 & 10.786 & \textbf{9.47} \\
\hlinewd{1.0pt}
\end{tabular}
\end{table}
\begin{table}[htbp]
\setlength\columnsep{0.1pt}
\centering
\caption{Speed-up of Twingo3D}
\label{speed-up_Twingo3D}
\begin{tabular}{lcccc}
\hlinewd{1.0pt}
problem & \#iter & CPU time (s) & GPU time (s) & speed-up \\
\hline
\multicolumn{3}{l}{\emph{P-BiCGSTAB}} & & \\
Twingo3D-0 & 563 & 1.85 & 1.008 & \textbf{1.84} \\
Twingo3D-1 & 1000 & 29.45 & 5.730 & \textbf{5.14} \\
Twingo3D-2 & 1000 & 295.66 & 37.670 & \textbf{7.85} \\
\hline
\multicolumn{3}{l}{\emph{P-BiCGSTAB(8)}} & & \\
Twingo3D-0 & 1000 & 31.2 & 20.970 & \textbf{1.49} \\
Twingo3D-1 & 1000 & 273.81 & 54.630 & \textbf{5.01} \\
Twingo3D-2 & 1000 & 2559.67 & 324.500 & \textbf{7.89} \\
\hline
\multicolumn{3}{l}{\emph{P-TFQMR}} & & \\
Twingo3D-0 & 366 & 1.34 & 0.626 & \textbf{2.14} \\
Twingo3D-1 & 954 & 30.4 & 5.438 & \textbf{5.59} \\
Twingo3D-2 & 1000 & 318.93 & 38.090 & \textbf{8.37} \\
\hlinewd{1.0pt}
\end{tabular}
\end{table}
The results corroborate the effectiveness of GPU compared to CPU for solving sparse linear systems. The speed-up grows when the size of the problems increase for all tests, i.e., for a finer mesh GPU is more effective compared to CPU.
For a finer mesh the assembled matrix turns into non appropriate size for memory of most of GPUs. To overcome this problem, one way consists in using domain decomposition method~\cite{SBG1996,quarteroni_domain_1999,toselli_domain_2004,NME:NME1620320604,magoules:journal-auth:21} based on iterative methods with interface conditions defined on the interface between the subdomains~\cite{magoules:journal-auth:16}.
The Schwarz method~\cite{magoulesf_contrib_3:Lions:1988:SAM,magoulesf_contrib_3:Lions:1989:SAM,magoulesf_contrib_3:Lions:1990:SAM,cai_overlapping_1998} is suitable for solving large size problem.
To accelerate the convergence, many references~\cite{magoulesf_contrib_3:Chevalier:1998:SMO,magoulesf_contrib_3:Gander:2000:OSM,magoules:journal-auth:23,magoules:journal-auth:18,magoules:journal-auth:14} show the importance of these interface conditions.
In order to implement this perspective for acoustic problems continuous optimized interface conditions between the subdomains must be implemented as in~\cite{magoules:journal-auth:6,magoules:journal-auth:10,magoules:journal-auth:9,magoules:journal-auth:13}.
Alternative discrete optimization techniques as introduced in~\cite{magoules:journal-auth:8,magoules:proceedings-auth:6,magoules:journal-auth:29,magoules:journal-auth:12,magoules:journal-auth:30,magoules:journal-auth:20,magoules:journal-auth:17} allow a fast and robust convergence of the Schwarz algorithm too.
In~\cite{cheikahamed:2013:inproceedings-4,ahamed2013stochastic}, the authors describe how domain decomposition method is effectively implemented on GPU and proved the robustness of Schwarz methods on a cluster of GPUs.
The extension to the complex number arithmetics double precision, of the iterative Krylov methods, to solve the local subproblems defined in each subdomains, leads to similar speed-up. For the Audi car compartment, a speed-up up to 9.2x is obtained for eight subdomains.

\section{Conclusion}
\label{sec:concluding_remarks}
In this paper we give an analysis of linear algebra operations together with their uses within iterative Krylov methods for solving acoustic problems on Graphics Processing Unit (GPU) with complex number arithmetics with double precision. Numerical tests have been carried out on two different system of accelerated generations of NVIDIA graphics card: GTX570 and Tesla K20c. A set of industrial matrices coming from the finite element discretization of acoustic problems modeled by the Helmholtz equation inside a car compartment are used to demonstrate the interest of using GPU device to perform linear algebra operations, and outline the robustness, performance and effectiveness of the proposed implementation.

\bibliography{bib/dcabes2014_acoustic_ac_fm,bib/MAGOULES-JOURNAL1,bib/MAGOULES-PROCEEDINGS1}

\bibliographystyle{abbrv}

\end{document}